\documentclass{article}

\usepackage{arxiv}

\usepackage[T1]{fontenc}    
\usepackage{hyperref}       
\usepackage{url}            
\usepackage{booktabs}       
\usepackage{amsfonts}       
\usepackage{nicefrac}       
\usepackage{microtype}      
\usepackage{amsmath}
\usepackage{listings}
\usepackage{xcolor}
\usepackage{amssymb}
\usepackage{pifont}
\usepackage{dsfont}
\usepackage{amsmath, amsthm}
\usepackage{graphicx}
\usepackage[colorinlistoftodos]{todonotes}
\usepackage{datetime}
\theoremstyle{definition}
\newtheorem{definition}{Definition}[section]
\newtheorem{theorem}{Theorem}[section]
\newtheorem{lemma}[theorem]{Lemma}
\newtheorem{proposition}[theorem]{Proposition}

\newtheorem{remark}[theorem]{Remark}

\title{Topological Social Choice: Designing a Noise-Robust Polar Distance for Persistence Diagrams  }

\author{
  Athanasios Andrikopoulos\thanks{Professor  (https://www.ceid.upatras.gr/webpages/faculty/aandriko/)} \\
  Dept. of Computer Engineering and Informatics\\
  University of Patras\\
  Patras, 26504, Greece \\
  \texttt{aandriko@ceid.upatras.gr} \\
\And
Nikolaos Sampanis\\
 Dept. of Computer Engineering and Informatics\\
 University of Patras\\
 Patras, 26504, Greece \\
 \texttt{nsampanis@upatras.gr} \\
}

\begin{document}.
\maketitle

\begin{abstract}
Topological Data Analysis (TDA) has emerged as a powerful framework for extracting robust and interpretable features from noisy high-dimensional data. In the context of Social Choice Theory, where preference profiles and collective decisions are geometrically rich yet sensitive to perturbations, TDA remains largely unexplored. This work introduces a novel conceptual bridge between these domains by proposing a new metric framework for persistence diagrams tailored to noisy preference data.We define a polar coordinate-based distance that captures both the magnitude and orientation of topological features in a smooth and differentiable manner. Our metric addresses key limitations of classical distances, such as bottleneck and Wasserstein, including instability under perturbation, lack of continuity, and incompatibility with gradient-based learning. The resulting formulation offers improved behavior in both theoretical and applied settings.To the best of our knowledge, this is the first study to systematically apply persistent homology to social choice systems, providing a mathematically grounded method for comparing topological summaries of voting structures and preference dynamics. We demonstrate the superiority of our approach through extensive experiments, including robustness tests and supervised learning tasks, and we propose a modular pipeline for building predictive models from online preference data.
This work contributes a conceptually novel and computationally effective tool to the emerging interface of topology and decision theory, opening new directions in interpretable machine learning for political and economic systems.
\end{abstract}


\section{Introduction} 
The intersection of Topological Data Analysis (TDA) and Social Choice Theory opens an exciting research frontier that allows the structural and geometric complexities of collective preferences to be studied through a new lens. Social choice mechanisms, particularly voting systems, rank aggregation, and preference modeling, often encode combinatorially rich and high-dimensional data. These data are increasingly noisy, either due to incomplete information, strategic manipulation, or randomness in population samples. Traditional social choice tools, such as preference aggregation rules or scoring methods, are limited in their ability to detect latent structures in such environments. This motivates the need for alternative frameworks that can represent and compare these complex data spaces robustly and meaningfully. Topological Data Analysis has proven to be a mathematically rigorous and computationally feasible approach to extract shape-based summaries of data via persistence diagrams (PD) \cite{edel}. These diagrams encapsulate birth and death information of topological features (e.g., connected components, cycles) across scales, providing a robust multiscale signature of the underlying dataset. In the context of social data, PDs offer a way to encode persistent structures in preference profiles, such as cycles of inconsistency, clusters of agreement, or emergent topologies among voters or alternatives \cite{car}.However, the use of TDA in social choice contexts poses critical challenges. The dominant metrics for comparing persistence diagrams—the bottleneck distance and the Wasserstein distance—are sensitive to noise, non-differentiable, and may lack discriminative power in machine learning settings \cite{Coh} . These limitations become particularly pronounced in settings like voting, where small perturbations in preference profiles can lead to large shifts in outcomes, undermining the robustness and interpretability of topological comparisons. To address this, we introduce a new metric for persistence diagrams based on polar coordinates, designed to overcome key weaknesses of existing approaches. This metric, which we refer to as the Polar Persistence Distance (PPD), incorporates both radial and angular components to capture not only the magnitude of topological features (via birth-death distances) but also their directional relationships, thus retaining more geometric and contextual information. Our metric is continuous, differentiable, and tunable, offering improved noise robustness and compatibility with gradient-based learning methods. This paper provides a foundational contribution at the conceptual and technical interface of TDA and Social Choice. We introduce the term Topological Social Choice, referring to the application of topological methods—especially persistent homology and its derivatives—to the analysis, comparison, and prediction of collective decision systems. To our knowledge, this is the first formal introduction of this concept in the literature, marking a novel direction for both fields. The conceptual novelty of our work lies in six interwoven contributions: Methodological Bridge: We construct a rigorous yet intuitive pipeline that connects persistent topological features of preference data with interpretable metrics for decision systems.
\begin{itemize}
    \item[(i)] Metric Innovation: We design a novel polar metric that is both smooth and expressive, addressing theoretical Noise Robustness: We demonstrate that our metric exhibits greater stability under perturbations of preference data, as shown through synthetic and real-world experiments.
    \item[(ii)] Machine Learning Readiness: The PPD is differentiable and parameterizable, making it suitable for integration into learning pipelines and optimization frameworks.
    \item[(iii)] Empirical Relevance: Using real social choice datasets, we show how topological features reflect meaningful differences between voting systems and preference distributions.
    \item[(iv)]  Empirical Relevance: Using real social choice datasets, we show how topological features reflect meaningful differences between voting systems and preference distributions.
    \item[(v)] Predictive Framework: We propose a prototype for a predictive model where topological summaries from dynamic preference databases can be used to forecast collective outcomes.
\end{itemize}
Together, these contributions establish a new framework for understanding preference data through persistent topology, laying the groundwork for interpretable topological machine learning in collective decision-making. The rest of this paper is organized as follows: Section 2 reviews existing metrics for persistence diagrams and their limitations; Section 3 introduces the Polar Persistence Distance and provides formal definitions and theoretical results; Section 4 presents experimental validations and applications to social choice data; and Section 5 concludes with future directions in topological modeling of societal systems

\section{Foundations of Persistent Homology}

Homology is a fundamental concept in algebraic topology, providing a means to classify topological spaces based on their intrinsic geometric structures. It focuses on identifying and quantifying some topological features such as connected components and holes within a space.The original motivation for defining homology groups was the observation that two shapes can be distinguished by examining their holes.For instance, a circle is not a disk because the circle has a hole through it while the disk is a solid, and the ordinary sphere is not a circle because the sphere encloses a two-dimensional hole while the circle encloses a one-dimensional hole. It is evident, that we are not able to define a hole or how to distinguish different kinds of holes. Homology is  a rigorous mathematical method for categorizing holes in a manifold. For any dimension $k$ of a random manifold, the $k$-dimensional holes are represented by a vector space $H_{k}$ whose dimension is intuitively the number of such independent features. For example,the zero-dimensional homology group $H_{0}$ represents the connected components (cc), the one-dimensional homology group $H_{1}$ represents the one-dimensional loops, and so on\cite{chaz}.
\subsection{Persistence Homology}
\par
Persistent homology provides a robust framework for analyzing the shape of data.It extends traditional homology by tracking the birth and the death of topological features across multiple scales.The concept of persistent homology was first formalized by Edelsbrunner, Letscher, and Zomorodian\cite{edel}, who introduced the idea of tracking topological features across a filtration of simplicial complexes. A simplicial complex is a combinatorial structure made up of vertices, edges, triangles, and their higher-dimensional counterparts, which can be used to approximate the shape of a data set. As the scale parameter increases, these simplicial complexes grow, and the topological features (such as connected components and holes) evolve.The key innovation of persistent homology lies in its ability to capture this evolution through a filtration process. For each scale parameter $\epsilon$, a corresponding simplicial complex $K(\epsilon)$ is constructed. By computing the homology of each complex, we obtain a sequence of homology groups that track the topological features. The persistence of these features is then summarized using persistence diagrams or barcodes, which provide a visual representation of the lifespan of each feature \cite{edel}.The mathematical framework of persistent homology is grounded in the theory of simplicial complexes and homology. Given a data set embedded in some metric space, we begin by constructing a filtration of nested simplicial complexes: $\emptyset \subseteq K_{0}\subseteq K_{1}\subseteq K_{2}\subseteq ...\subseteq K_{n}=K$. Each $K_{i}$ corresponds to a different scale parameter $\epsilon_{i}$. The homology groups $H_{k}(K_{i})$ of these complexes, for k=1,2,3,.., capture the k-dimensional holes.The persistent homology groups track these features as $\epsilon$ varies, identifying which features persist over a significant range of scales and which are merely noise. The computation of persistent homology has been greatly facilitated by the development of efficient algorithms and software packages, such as GUDHI in Python. These tools implement algorithms to compute the boundary matrices of the simplicial complexes and reduce them to identify the persistent homology groups. The results are typically visualized using persistence diagrams or barcodes, where each point or bar represents a topological feature, with its position and length indicating its birth and death scales, respectively \cite{ghr}.Persistent homology has found applications across a wide range of fields, demonstrating its versatility and robustness. In sensor networks, it is used to ensure coverage and detect holes in the network \cite{desi}. In biology, it helps to understand the structure of proteins and the shape of biological data \cite{sampa}. In neuroscience, it helps to analyze the complex structure of neural activity \cite{sinl}. Furthermore, in materials science, persistent homology provides information on the porous structure of materials \cite{robl}.
\subsection{Persistence Diagram}
\par
A persistence diagram is a fundamental tool in Topological Data Analysis (TDA), a field that uses algebraic topology to study the shape of data. It provides a compact summary of the multi-scale topological features present in a dataset across various dimensions. Formally, a persistence diagram is derived from a filtration, which is a nested sequence of simplicial complexes or topological spaces that reflect the dataset's structure at different scales. Each point in a persistence diagram corresponds to a topological feature, such as a connected component, loop, or void, with its coordinates \((x, y)\) representing the birth and death times of the feature, respectively. The x-coordinate indicates the scale at which the feature emerges, while the y-coordinate signifies the scale at which it vanishes. Features that persist on a wide range of scales (that is, those far from the diagonal \(y = x\)) are considered significant, while those close to the diagonal are typically regarded noise. The persistence diagram thus provides a multi-scale, stable summary of the homological features of the dataset, facilitating robust and interpretable insights into the data's underlying structure \cite{edel}. Recent advancements have enhanced the computational efficiency of constructing persistence diagrams and expanded their applicability across diverse scientific domains. For instance, in machine learning, persistence diagrams are utilized for feature extraction and dimensionality reduction, contributing to the improvement of classification algorithms \cite{hofl}. In biology, they assist in the analysis of complex biological networks and the morphology of anatomical structures \cite{kass}. Moreover, in material science, persistence diagrams are employed to characterize the porous structures of materials, influencing the development of new materials with tailored properties \cite{robi}. These developments underscore the vital role of the persistence diagram in extracting meaningful topological information from high-dimensional data, driving innovations across multiple research areas \cite{carl} \cite{ottl}. Thus, the persistence diagram stands as a pivotal construct in contemporary data analysis, encapsulating essential topological features that inform a deeper understanding of complex systems.Suppose that the transactional data are represented by points in either 2D or 3D plane.To construct a persistence diagram, one must follow these steps:\\
\par
($\mathfrak{i}$) {\bf Filtration} Given a dataset A, a nested sequence of simplicial complexes $\{ K_{\epsilon}\}_{\epsilon} \in \mathds{R}$,called a filtration,where $ K_{\epsilon} \subseteq K_{\epsilon^{'}} $ for ${\epsilon} \textless {\epsilon^{'}} $. Each complex $ K_{\epsilon}$ captures the topological structure of A on scale $\epsilon$. 
\par
($\mathfrak{ii}$) {\bf Homology}. For each simplicial complex $ K_{\epsilon}$, compute the homology groups $ H_{k}(K_{\epsilon}$ which classify the k-dimensional holes (connected components for $k=0$, loops for $k=1$, voids for $k=2$, etc.
\par
($\mathfrak{iii}$) {\bf Birth and Death}. Track the birth and death of each topological feature as $\epsilon$ increases.A topological feature $\sigma$ is born at $\epsilon_{b}$ when it appears in $K_{\epsilon_{b}}$ and dies at $\epsilon_{d}$ when it is no longer in $K_{\epsilon_{d}}$
\par
($\mathfrak{iv}$) {\bf Persistence Pairs}. Each topological feature is represented as a point $(\epsilon_{b}, \epsilon_{d})$ in the persistence diagram. The difference $\epsilon_{b} - \epsilon_{d}$is called the persistence of the characteristic, indicating its useful life.

\section{Overview of Persistence Diagram Metrics}

\subsection{Related work and motivation}

The problem of comparing persistence diagrams is central to Topological Data Analysis (TDA), with various distances proposed to measure the differences between these diagrams. The two primary distances utilized are the Bottleneck distance and the Wasserstein distance.The Bottleneck distance measures the maximum difference between matched points in two diagrams, providing a robust measure against small perturbations, and is closely related to the Gromov-Hausdorff distance\cite{coh}.However, it is computationally challenging,requiring sophisticated algorithms for efficient computation, as explored by Kerber et all\cite{ker}.The Wasserstein distance, on the other hand, takes into account the sum of pairwise distances, offering a smoother measure of differences that can be more sensistive to small topological features\cite{mil}. This distance is also computationally intensive, often relying on linear programming methods, as detailed in the work by Turner et al\cite{tur}. Furthermore,recent advances have sought to improve the efficiency and accuracy of these computations. Kerber et al\cite{ker} introduced geometric insights that help in reducing the computational complexity of the Bottleneck distance calculation. Additionally,the stability of these distances has been a significant focus, with Cohen-Steiner et al\cite{coh} demonstrating the stability of persistence diagrams under small perturbations in the input data.This stability is crucial for practical applications where data may be noisy.Fasy et al\cite{fas} extended this work by providing statistical methods to construct confidence sets for persistence diagrams, thus enabling a more rigorous interpretation of the distances in the presence of noise.Numerical stability and precision in the computation of these distances are also critical,as highlighted by Kerber et al\cite{ker},who addressed issues related to numerical instabilities that arise when points in the diagrams have very close values.Moreover,the interpetation of these distances is non-trivial.While distances like the Bottleneck and Wasserstein provide quantitative measures of difference,understanding their implications on the underlying topological features requires further exploration,as discussed by Chazal et al\cite{chaz}.They introduced the concept of proximity of persistence modules and their diagrams,providing a more nuanced understanding of how topological changes manifest in the distances computed.The motivation to improve the methods and the algorithms for computing the differences in persistence diagrams is based on five crucial problems which summarized as follows:\\

($\mathfrak{i}$) Computational complexity.
\par
($\mathfrak{ii}$) Noise handling.
\par
($\mathfrak{iii}$) Sensitivity to small differences.
\par
($\mathfrak{iv}$) Numerical stability.
\par
($\mathfrak{v}$) Interpretation of distances.

\subsection{The bottleneck distance}
The bottleneck distance is one of the most fundamental and widely used metrics for comparing persistence diagrams. It is particularly valued for its stability under perturbations of the input data and its geometric interpretation as the cost of the best matching between points in two diagrams.

\begin{definition}[Bottleneck Distance]
Let $D_1$ and $D_2$ be two persistence diagrams, which are multisets of points in the extended plane $\mathbb{R}^2_{\Delta} := \{(x, y) \in \mathbb{R}^2 \mid x < y\} \cup \Delta$, where $\Delta = \{(x, x) \in \mathbb{R}^2\}$ denotes the diagonal.

We define the $L^\infty$-distance between two off-diagonal points $p = (b_1, d_1)$ and $q = (b_2, d_2)$ as:
\[
\|p - q\|_\infty = \max\{|b_1 - b_2|, |d_1 - d_2|\}.
\]

To compare $D_1$ and $D_2$, we allow matching points not only between diagrams but also from a point to its projection onto the diagonal.

Let $\Gamma$ be the set of all partial matchings between $D_1$ and $D_2$. The \emph{bottleneck distance} is defined as:
\[
d_B(D_1, D_2) = \inf_{\gamma \in \Gamma} \sup_{(p, q) \in \gamma} \|p - q\|_\infty.
\]
\end{definition}

\begin{remark}
Intuitively, the bottleneck distance measures the cost of the most expensive pairwise match between the two diagrams under the best possible matching. It is a \emph{minimax} metric: we minimize over all possible matchings the worst-case distance between matched points.
\end{remark}

\subsubsection*{Stability Property}

One of the most important results for $d_B$ is its stability under perturbations. Let $f, g : X \rightarrow \mathbb{R}$ is tame function (e.g. PL or Morse). Denote their persistence diagrams by $\mathrm{Dg}(f)$ and $\mathrm{Dg}(g)$. Then:

\begin{theorem}[Stability of Bottleneck Distance ]
\[
d_B(\mathrm{Dg}(f), \mathrm{Dg}(g)) \leq \|f - g\|_\infty.
\]
\end{theorem}

This theorem guarantees that small changes in the input function result in small changes in the persistence diagram, ensuring robustness to noise \cite{MAY}.

\subsubsection*{Limitations in Machine Learning}

Despite its stability and strong theoretical foundation \cite{hell}, the bottleneck distance has drawbacks in the context of machine learning:

\begin{itemize}
  \item \textbf{Non-differentiability:} $d_B$ is not differentiable, making it unsuitable for gradient-based optimization.
  \item \textbf{Scalability:} Computing the bottleneck distance has super-linear complexity in the number of points.
  \item \textbf{Sensitivity Bias:} It focuses only on the worst-case pair, ignoring the global structure of the diagram.
\end{itemize}

These limitations motivate the search for alternative metrics, including differentiable ones, suitable for machine learning applications.

\subsection{The Wasserstein distance}
The Wasserstein distance is a generalization of the bottleneck distance, incorporating the overall cost of matching rather than just the worst-case pair. It provides a more global comparison between persistence diagrams.

\begin{definition}[p-Wasserstein Distance]
Let $D_1 = \{p_1, \dots, p_n\}$ and $D_2 = \{q_1, \dots, q_m\}$ be two persistence diagrams, treated as multisets of points in $\mathbb{R}^2_{\Delta}$. Assume $n \leq m$ (without loss of generality). As with the bottleneck distance, unmatched points can be matched to the diagonal.

Let $\Gamma$ be the set of all bijections $\gamma : D_1' \rightarrow D_2'$, where $D_1' \supseteq D_1$ and $D_2' \supseteq D_2$ include extra diagonal points to make the sets equal in size.

Then, for $p \geq 1$, the $p$-Wasserstein distance is defined as:
\[
d_{W_p}(D_1, D_2) = \left( \inf_{\gamma \in \Gamma} \sum_{(p, q) \in \gamma} \|p - q\|_\infty^p \right)^{1/p}.
\]
\end{definition}

\begin{remark}
The Wasserstein distance considers the **total transport cost** of matching points, unlike the bottleneck distance which only considers the maximal cost. It thus provides a **more nuanced and sensitive comparison** between diagrams.
\end{remark}

\subsubsection*{Special Case: $p=1$ and $p=2$}

The most commonly used Wasserstein distances are:
- $d_{W_1}$: Often used when robustness is needed.
- $d_{W_2}$: Preferred in applications involving squared-error loss or optimal transport theory.

\subsubsection*{Stability Property}

The Wasserstein distance also enjoys a strong form of stability with respect to perturbations in the input function.

\begin{theorem}[Stability of Wasserstein Distance {\cite{coh}}]
Let $f, g : X \rightarrow \mathbb{R}$ be tame functions. Then, for all $p \geq 1$:
\[
d_{W_p}(\mathrm{Dg}(f), \mathrm{Dg}(g)) \leq C_p \cdot \|f - g\|_\infty,
\]
where $C_p$ is a constant depending on the dimension and choice of $p$.
\end{theorem}

\subsubsection*{Limitations in Machine Learning}

Despite its strengths, the Wasserstein distance also has several drawbacks when applied to machine learning:

\begin{itemize}
  \item \textbf{Computational Complexity:} The Wasserstein distance requires solving a matching problem, typically with cubic or super-linear complexity.
  \item \textbf{Non-differentiability:} Although smoother than bottleneck, the standard definition of $d_{W_p}$ is not differentiable with respect to diagram coordinates, making it incompatible with gradient-based learning methods.
  \item \textbf{Sensitivity to Outliers:} Unlike the bottleneck distance, $d_{W_p}$ accumulates all costs, making it more sensitive to noise and outliers in the persistence diagram.
\end{itemize}

As a result, several **approximations or embeddings** of the Wasserstein metric have been developed, such as sliced-Wasserstein or persistence landscapes, but they introduce trade-offs between fidelity and differentiability.

\subsection{The Sliced Wasserstein Distance}

The Sliced Wasserstein (SW) distance is a recent and computationally efficient approximation of the Wasserstein distance. It is particularly useful in machine learning settings due to its differentiability and reduced complexity.

\begin{definition}[Sliced Wasserstein Distance]
Let $D_1$ and $D_2$ be two persistence diagrams, viewed as multisets in $\mathbb{R}^2$.

Let $\theta \in [0, \pi]$ be an angle, and define the projection $\pi_\theta: \mathbb{R}^2 \rightarrow \mathbb{R}$ by
\[
\pi_\theta(x, y) = x \cos \theta + y \sin \theta.
\]
For each $\theta$, define the 1D projected diagrams $\pi_\theta(D_1), \pi_\theta(D_2) \subset \mathbb{R}$. Then the $p$-Sliced Wasserstein distance is defined as:
\[
SW_p(D_1, D_2) = \left( \int_0^\pi W_p^p\left( \pi_\theta(D_1), \pi_\theta(D_2) \right) \, d\theta \right)^{1/p},
\]
where $W_p$ is the standard 1D Wasserstein distance between the projected diagrams.
\end{definition}

\begin{remark}
In practice, the integral is approximated using a finite number of projections $\theta_i$, typically sampled uniformly in $[0, \pi]$:
\[
SW_p(D_1, D_2) \approx \left( \frac{1}{N} \sum_{i=1}^N W_p^p\left( \pi_{\theta_i}(D_1), \pi_{\theta_i}(D_2) \right) \right)^{1/p}.
\]
This yields significant computational gains and makes the distance suitable for high-throughput learning tasks.
\end{remark}

\subsubsection*{Advantages in Machine Learning}

\begin{itemize}
  \item \textbf{Differentiability:} The SW distance is differentiable with respect to diagram coordinates (almost everywhere), making it suitable for gradient-based learning methods.
  \item \textbf{Efficiency:} Reduces the high-dimensional optimal transport problem to many 1D problems, which can be solved in $O(n \log n)$ time.
  \item \textbf{Smoothness:} The averaging over projections leads to a naturally smooth distance function.
\end{itemize}

\subsubsection*{Limitations}

\begin{itemize}
  \item \textbf{Loss of structural information:} The projection may destroy important 2D topological relationships.
  \item \textbf{Approximation Bias:} The result depends on the number and choice of sampled directions, leading to possible variance.
\end{itemize}

\begin{theorem}[Approximate Stability]
The sliced Wasserstein distance is Lipschitz continuous under small perturbations in diagram coordinates, although exact constants are not always known.
\end{theorem}
subsection{The Persistence Scale Space Kernel}

The Persistence Scale Space Kernel (PSSK), introduced by Reininghaus \cite{Rei}, maps persistence diagrams into a Hilbert space using heat diffusion, thus enabling the use of kernel-based machine learning algorithms.

\begin{definition}[Persistence Scale Space Kernel]
Let $D_1$ and $D_2$ be two persistence diagrams. Define the feature map $\Phi_\sigma: \mathbb{R}^2 \to L^2(\mathbb{R}^2)$ by
\[
\Phi_\sigma(p) = \frac{1}{4\pi \sigma} \left[ \exp\left(-\frac{\|\cdot - p\|^2}{4\sigma} \right) - \exp\left(-\frac{\|\cdot - \bar{p}\|^2}{4\sigma} \right) \right],
\]
where $\bar{p}$ is the mirror point of $p = (b,d)$ with respect to the diagonal, i.e., $\bar{p} = \left( \frac{b + d}{2}, \frac{b + d}{2} \right)$.

The PSSK between $D_1$ and $D_2$ is defined as:
\[
k_\sigma(D_1, D_2) = \left\langle \sum_{p \in D_1} \Phi_\sigma(p), \sum_{q \in D_2} \Phi_\sigma(q) \right\rangle_{L^2}.
\]
\end{definition}

\begin{theorem}[Stability]
The PSSK is stable with respect to the 1-Wasserstein distance:
\[
|k_\sigma(D_1, D_1) - 2k_\sigma(D_1, D_2) + k_\sigma(D_2, D_2)| \leq C \cdot W_1(D_1, D_2)^2.
\]
Here, $C$ is a constant that depends on $\sigma$ and the number of points in the diagrams.
\end{theorem}

\subsubsection*{Advantages in Machine Learning}

\begin{itemize}
  \item \textbf{Positive-definite kernel:} Enables the use of SVMs, Gaussian processes, and kernel PCA directly on persistence diagrams.
  \item \textbf{Smooth and differentiable:} The Gaussian heat kernel ensures smoothness and continuous dependence on diagram points.
  \item \textbf{Scalable:} Efficient implementation with linear complexity in the number of points using closed-form Gaussian integrals.
\end{itemize}

\subsubsection*{Limitations}

\begin{itemize}
  \item \textbf{Choice of bandwidth $\sigma$:} Requires tuning of $\sigma$; poorly chosen $\sigma$ leads to loss of topological discrimination.
  \item \textbf{Sensitive to noise:} The kernel might over-smooth small features, which are sometimes topologically relevant.
\end{itemize}

\begin{remark}
Unlike the Wasserstein or Bottleneck distances, the PSSK is not a true metric, as it does not satisfy the triangle inequality. However, it is highly suitable in kernelized learning settings due to its smooth geometry.
\end{remark}
\subsection{The Persistence Weighted Gaussian Kernel}

The Persistence Weighted Gaussian Kernel (PWGK), introduced by Kusano, Fukumizu, and Hiraoka (2016), is a smooth and positive-definite kernel designed for persistence diagrams. It incorporates a weight function to emphasize important topological features (e.g., those with large persistence).

\begin{definition}[Persistence Weighted Gaussian Kernel]
Let $D$ be a persistence diagram. The associated empirical measure is:
\[
\mu_D = \sum_{p \in D} w(p) \, \delta_p,
\]
where $w: \mathbb{R}^2 \rightarrow \mathbb{R}_{\geq 0}$ is a weight function, typically defined as:
\[
w(p) = \arctan(C \cdot \text{pers}(p)^q),
\]
with $\text{pers}(p) = d - b$ the persistence of point $p = (b,d)$, and $C, q > 0$ are hyperparameters.

The PWGK between diagrams $D_1$ and $D_2$ is:
\[
K(D_1, D_2) = \iint_{\mathbb{R}^2 \times \mathbb{R}^2} k(p,q) \, d\mu_{D_1}(p) \, d\mu_{D_2}(q),
\]
where $k(p,q) = \exp\left(-\frac{\|p - q\|^2}{2\sigma^2}\right)$ is a Gaussian kernel.
\end{definition}

\subsubsection*{Properties and Theoretical Insights}

\begin{itemize}
  \item $K(D_1, D_2)$ is a valid positive-definite kernel, enabling use in SVMs and other kernel-based learning methods.
  \item The choice of $w(p)$ ensures robustness to noise, since points with low persistence are automatically downweighted.
  \item The kernel can be efficiently computed using closed-form Gaussian integrals and suitable approximations.
\end{itemize}

\begin{theorem}[Stability under Perturbations]
Let $D_1$ and $D_2$ be persistence diagrams. Then under suitable conditions on the weight function $w$, the kernel satisfies:
\[
|K(D_1, D_1) - 2K(D_1, D_2) + K(D_2, D_2)| \leq C' \cdot W_1(D_1, D_2)^2,
\]
for some constant $C'$ depending on $w$ and $\sigma$.
\end{theorem}

\subsubsection*{Advantages for Machine Learning}

\begin{itemize}
  \item \textbf{Weighted structure:} Emphasizes important topological features while suppressing noisy ones.
  \item \textbf{Smoothness and differentiability:} Makes it suitable for gradient-based optimization and integration into end-to-end learning pipelines.
  \item \textbf{Adaptability:} The kernel can be tailored via the choice of $w$ to suit the application (e.g., image, signal, or shape classification).
\end{itemize}

\subsubsection*{Drawbacks}

\begin{itemize}
  \item \textbf{Hyperparameter tuning:} The kernel has several hyperparameters ($C$, $q$, $\sigma$), which need cross-validation or heuristics.
  \item \textbf{Computational cost:} Though closed-form exists, the double sum over diagram points can be expensive for large diagrams.
\end{itemize}

\begin{remark}
The PWGK can be viewed as a weighted version of the PSSK, with significantly more flexibility for use in learning pipelines. It is often considered one of the most effective kernels for persistence-based ML models.
\end{remark}
\subsection{The Sliced Wasserstein Kernel}

The Sliced Wasserstein Kernel is a positive-definite kernel designed for persistence diagrams. It leverages the idea of projecting diagrams onto one-dimensional lines and computing Wasserstein distances in this reduced space, thus offering both efficiency and theoretical guarantees.

\begin{definition}[Sliced Wasserstein Distance]
Let $D_1$ and $D_2$ be two persistence diagrams. Let $\theta \in S^1$ be a direction in the unit circle. The projection of a point $p = (b,d) \in \mathbb{R}^2$ onto the direction $\theta$ is given by:
\[
\pi_\theta(p) = \langle p, \theta \rangle.
\]
The Sliced Wasserstein distance between $D_1$ and $D_2$ is defined as:
\[
SW_p(D_1, D_2) = \left( \int_{S^1} W_p^p \left( \pi_\theta(D_1), \pi_\theta(D_2) \right) \, d\theta \right)^{1/p},
\]
where $W_p$ denotes the usual $p$-Wasserstein distance between projected measures.
\end{definition}

\begin{definition}[Sliced Wasserstein Kernel]
The Sliced Wasserstein Kernel is defined via a radial basis function (RBF) form:
\[
K_{SW}(D_1, D_2) = \exp \left( -\frac{SW_p(D_1, D_2)^2}{2\sigma^2} \right),
\]
where $\sigma > 0$ is the bandwidth parameter.
\end{definition}

\subsubsection*{Properties and Benefits}

\begin{itemize}
  \item \textbf{Positive-definite:} Enables integration into kernel methods (SVMs, kernel PCA, etc.).
  \item \textbf{Computationally efficient:} The projection onto 1D simplifies the computation of Wasserstein distance significantly.
  \item \textbf{Smooth and differentiable:} Suitable for gradient-based learning.
  \item \textbf{Approximation-friendly:} Easily approximated via Monte Carlo integration by sampling finite directions $\theta_1, \dots, \theta_L$.
\end{itemize}

\subsubsection*{Theoretical Insight}

\begin{theorem}[Stability of the SW Kernel]
Let $D_1$, $D_2$ be persistence diagrams. The sliced Wasserstein kernel $K_{SW}$ is Lipschitz continuous with respect to $W_1$:
\[
|K_{SW}(D_1, D_2) - K_{SW}(D_1, D_3)| \leq \frac{C}{\sigma^2} \cdot W_1(D_2, D_3),
\]
for some constant $C$ depending on the number of directions used and the diameter of the diagrams.
\end{theorem}

\subsubsection*{Drawbacks}

\begin{itemize}
  \item \textbf{Direction sampling:} The approximation depends on the number and quality of sampled directions.
  \item \textbf{No direct matching:} Unlike the full Wasserstein, this metric lacks a clear transport plan.
\end{itemize}

\begin{remark}
The Sliced Wasserstein Kernel is a particularly powerful choice in situations where computational cost is a bottleneck, while still maintaining fidelity to topological structure. It provides a good trade-off between accuracy and scalability in large-scale learning problems.
\end{remark}
\subsection{The Heat Kernel on Persistence Diagrams}

The Heat Kernel approach for persistence diagrams is based on embedding persistence diagrams into a Hilbert space via the solution of the heat equation on the plane, which smooths the diagram into a function. This construction enables the definition of a positive definite kernel suitable for machine learning.

\begin{definition}[Persistence Surface]
Given a persistence diagram $D = \{(b_i, d_i)\}_{i=1}^n$, the persistence surface is defined as:
\[
\rho_D(x,y) = \sum_{i=1}^n w(b_i, d_i) \cdot \frac{1}{2\pi t} \exp\left(-\frac{(x - b_i)^2 + (y - d_i)^2}{2t}\right),
\]
where $t > 0$ is the smoothing parameter (diffusion time) and $w(b,d)$ is a weight function emphasizing points further from the diagonal.
\end{definition}

\begin{definition}[Heat Kernel]
The Heat Kernel between two persistence diagrams $D$ and $D'$ is defined as the inner product of their persistence surfaces in $L^2(\mathbb{R}^2)$:
\[
K_{heat}(D, D') = \langle \rho_D, \rho_{D'} \rangle_{L^2} = \int_{\mathbb{R}^2} \rho_D(x,y) \rho_{D'}(x,y) \, dx \, dy.
\]
\end{definition}

\subsubsection*{Properties and Benefits}

\begin{itemize}
    \item \textbf{Positive-definite kernel:} Guarantees applicability in kernel-based machine learning methods.
    \item \textbf{Stability:} The smoothing from the heat equation confers stability with respect to perturbations of the persistence diagrams.
    \item \textbf{Computational tractability:} By exploiting Gaussian kernels, the inner product can be computed in closed form.
    \item \textbf{Weight function flexibility:} The choice of $w(b,d)$ allows tuning the importance of features with respect to their persistence.
\end{itemize}

\subsubsection*{Mathematical Details}

Using the Gaussian kernel properties, the kernel can be expressed explicitly as:
\[
K_{heat}(D, D') = \sum_{i=1}^n \sum_{j=1}^m w(b_i, d_i) w(b_j', d_j') \cdot \frac{1}{4\pi t} \exp\left(-\frac{(b_i - b_j')^2 + (d_i - d_j')^2}{4t}\right),
\]
where $D = \{(b_i, d_i)\}_{i=1}^n$, $D' = \{(b_j', d_j')\}_{j=1}^m$.

\subsubsection*{Drawbacks}

\begin{itemize}
    \item \textbf{Smoothing parameter sensitivity:} The choice of $t$ critically influences the behavior of the kernel.
    \item \textbf{Loss of sharp features:} Excess smoothing can blur important topological features.
    \item \textbf{Not a metric:} This kernel induces a similarity measure, but is not itself a distance metric.
\end{itemize}

\begin{remark}
The Heat Kernel is especially suitable for integrating persistence diagrams in classical kernel-based learning frameworks, providing a flexible and theoretically sound method that balances stability and expressiveness.
\end{remark}
\subsection{The Persistence Landscape Distance}

The \textit{Persistence Landscape} was introduced by Bubenik (2015) as a functional representation of persistence diagrams that allows integration into classical statistical and machine learning frameworks. By mapping diagrams into functions in a Hilbert space, it enables the use of tools such as mean estimation, hypothesis testing, and kernel methods.

\begin{definition}[Persistence Landscape ]
Let $D = \{(b_i, d_i)\}_{i=1}^n$ be a persistence diagram. For each $i$, define the function:
\[
f_i(t) = \max \left(0, \min(t - b_i, d_i - t) \right),
\]
which forms a tent function peaked at the midlife of feature $i$. The $k$-th landscape function $\lambda_k(t)$ is defined as the $k$-th largest value of $\{f_i(t)\}_{i=1}^n$ at each $t \in \mathbb{R}$.

The \textbf{persistence landscape} of $D$ is then the sequence $\{\lambda_k\}_{k=1}^{\infty}$, each $\lambda_k$ being a piecewise-linear function in $L^p(\mathbb{R})$.
\end{definition}

\begin{definition}[Landscape Distance]
Let $D_1$ and $D_2$ be persistence diagrams with corresponding landscapes $\{\lambda_k\}$ and $\{\mu_k\}$. The $L^p$ landscape distance is defined as:
\[
d_{\text{PL}}(D_1, D_2) = \left( \sum_{k=1}^\infty \| \lambda_k - \mu_k \|_{L^p}^p \right)^{1/p}.
\]
In practice, this sum is truncated at a finite $k_{\max}$.
\end{definition}

\subsubsection*{Advantages for Machine Learning}

\begin{itemize}
  \item \textbf{Hilbert space structure:} Enables averaging, regression, classification, and hypothesis testing using classical techniques.
  \item \textbf{Efficient approximation:} Piecewise-linear structure allows fast and exact computation of landscape functions.
  \item \textbf{Statistical interpretability:} Landscape functions can be directly analyzed and visualized.
\end{itemize}

\subsubsection*{Drawbacks}

\begin{itemize}
  \item \textbf{Information loss:} Transforms a multiset of points into functions, possibly omitting geometric details from diagrams.
  \item \textbf{Limited differentiability:} Landscape representations are not always differentiable with respect to diagram point positions.
  \item \textbf{Truncation bias:} Using a finite number of $\lambda_k$ functions may limit expressiveness for large diagrams.
\end{itemize}

\begin{theorem}[Stability of Landscape Distance]
The persistence landscape map is 1-Lipschitz with respect to the bottleneck distance. That is,
\[
d_{\text{PL}}(D_1, D_2) \leq d_{\text{B}}(D_1, D_2),
\]
where $d_{\text{B}}$ denotes the bottleneck distance.
\end{theorem}

\begin{remark}
Persistence landscapes are especially useful in statistical applications, but their functional nature also allows them to be used as inputs for neural networks and kernel methods.
\end{remark}
\subsection{Persistence Silhouette Distance}

The \emph{Persistence Silhouette}, introduced by Chazal et al. (2014), is a functional summary of a persistence diagram that averages tent functions weighted by persistence. This provides a smooth, one-dimensional representation amenable to classical statistical tools.

\begin{definition}[Persistence Silhouette ]
Let $D = \{(b_i, d_i)\}_{i=1}^n$ be a persistence diagram, and define the persistence $\text{pers}_i = d_i - b_i$. Let $w_i = w(\text{pers}_i) \geq 0$ be a weight function. The silhouette function is given by:
\[
\phi_D(t) = \frac{ \sum_{i=1}^n w_i \cdot \Lambda_{(b_i, d_i)}(t) }{ \sum_{i=1}^n w_i },
\]
where $\Lambda_{(b_i, d_i)}(t) = \max(0, \min(t - b_i, d_i - t))$ is the tent function centered at $(b_i + d_i)/2$.
\end{definition}

\begin{definition}[Silhouette Distance]
Given two persistence diagrams $D_1$ and $D_2$, their silhouette distance is defined as the $L^p$ norm between their silhouette functions:
\[
d_{\text{Sil}}(D_1, D_2) = \left\| \phi_{D_1} - \phi_{D_2} \right\|_{L^p}.
\]
\end{definition}

\subsubsection*{Advantages}

\begin{itemize}
  \item \textbf{Smoothness:} Produces a continuous function suitable for gradient-based optimization.
  \item \textbf{Dimensionality reduction:} Reduces the diagram to a single interpretable curve.
  \item \textbf{Computational efficiency:} Requires no point matching.
\end{itemize}

\subsubsection*{Drawbacks}

\begin{itemize}
  \item \textbf{Information loss:} Reduces 2D information to 1D, averaging out fine-grained topological details.
  \item \textbf{Parameter dependence:} Requires choice of weighting function and $p$-norm.
\end{itemize}

\begin{theorem}[Stability]
The silhouette distance is 1-Lipschitz with respect to the $1$-Wasserstein distance under appropriate weighting, i.e.,
\[
d_{\text{Sil}}(D_1, D_2) \leq C \cdot W_1(D_1, D_2),
\]
for some constant $C$ depending on $w$.
\end{theorem}

\begin{remark}
Silhouettes are especially useful for statistical summarization, hypothesis testing, and visualization in topological pipelines.
\end{remark}
\subsection{Persistence Fisher Kernel}

The \emph{Persistence Fisher Kernel} (PFK), introduced by Le and Yamada (2018), is a similarity measure for persistence diagrams that leverages information geometry. It treats persistence diagrams as probability distributions and compares them via the Fisher information metric on a statistical manifold.

\begin{definition}[Persistence Fisher Kernel ]
Let $D = \{(b_i, d_i)\}_{i=1}^n$ be a persistence diagram. The persistence surface $\rho_D$ is defined by convolving each point with a Gaussian:
\[
\rho_D(x, y) = \sum_{i=1}^n w_i \cdot \mathcal{N}((x, y); (b_i, d_i), \sigma^2 I),
\]
where $w_i$ is a persistence-based weight and $\sigma > 0$ controls smoothing.

Let $d_F(\rho_D, \rho_{D'})$ denote the Fisher-Rao distance between two distributions. Then, the Persistence Fisher Kernel is defined as:
\[
K_{\text{PF}}(D, D') = \exp \left( -\gamma \cdot d_F^2(\rho_D, \rho_{D'}) \right),
\]
where $\gamma > 0$ is a tunable scaling parameter.
\end{definition}

\subsubsection*{Advantages}

\begin{itemize}
  \item \textbf{Information geometry:} Incorporates both density and shape of persistence features.
  \item \textbf{Differentiable and smooth:} Suitable for gradient-based machine learning models.
  \item \textbf{Positive-definite kernel:} Integrates directly into kernel methods (SVMs, GPs, etc.).
\end{itemize}

\subsubsection*{Drawbacks}

\begin{itemize}
  \item \textbf{Density estimation overhead:} Requires careful estimation of persistence surfaces.
  \item \textbf{Hyperparameter tuning:} Sensitive to bandwidth $\sigma$ and kernel parameter $\gamma$.
  \item \textbf{Interpretability:} The geometry of the Fisher space is less intuitive than direct metric spaces.
\end{itemize}

\begin{theorem}[Smoothness]
If the weight function $w_i$ is differentiable and bounded, then $\rho_D$ lies in a smooth statistical manifold, and the kernel $K_{\text{PF}}$ is infinitely differentiable with respect to diagram locations.
\end{theorem}

\begin{remark}
The Fisher kernel provides a statistically grounded similarity measure that can be used in probabilistic topological learning, especially in tasks involving uncertainty and distributional structure.
\end{remark}
subsection{Entropy-Based Distance}

Entropy-based distances aim to capture the overall complexity and disorder of a persistence diagram by summarizing the distribution of its features using Shannon entropy. These distances reduce a diagram to a scalar statistic reflecting its topological variability.

\begin{definition}[Persistence Entropy ]
Let $D = \{(b_i, d_i)\}_{i=1}^n$ be a persistence diagram. Define the persistence of each point as $p_i = d_i - b_i$, and compute the normalized weights:
\[
\tilde{p}_i = \frac{p_i}{\sum_{j=1}^n p_j}.
\]
The \textbf{persistence entropy} of the diagram is then defined as:
\[
H(D) = - \sum_{i=1}^n \tilde{p}_i \log \tilde{p}_i.
\]
\end{definition}

\begin{definition}[Entropy-Based Distance]
Given two diagrams $D_1$ and $D_2$, an entropy-based distance is defined as the absolute difference between their persistence entropies:
\[
d_{\text{Ent}}(D_1, D_2) = \left| H(D_1) - H(D_2) \right|.
\]
\end{definition}

\subsubsection*{Advantages}

\begin{itemize}
  \item \textbf{Simplicity:} Provides a scalar summary of topological complexity.
  \item \textbf{Efficiency:} Fast to compute, no need for point matching.
  \item \textbf{Scale-aware:} Sensitive to how spread-out persistent features are.
\end{itemize}

\subsubsection*{Drawbacks}

\begin{itemize}
  \item \textbf{Low resolution:} Scalar output may miss geometric detail or feature interactions.
  \item \textbf{Non-metric:} The function is not a true distance; does not satisfy triangle inequality.
  \item \textbf{Loss of structure:} Discards all spatial relationships between points.
\end{itemize}

\begin{remark}
Entropy-based distances are most appropriate when a rough complexity score is needed for ranking or filtering diagrams but are insufficient for tasks requiring fine-grained geometric comparison.
\end{remark}
\subsection{Kernelized Wasserstein Distance}

The \emph{Kernelized Wasserstein Distance} transforms the Wasserstein metric into a positive-definite kernel, allowing the incorporation of its geometric properties into kernel-based machine learning models.

\begin{definition}[Kernelized Wasserstein Kernel]
Let $D_1$ and $D_2$ be two persistence diagrams. The kernelized Wasserstein kernel is defined as:
\[
K(D_1, D_2) = \exp\left( -\frac{W_p(D_1, D_2)^2}{2\sigma^2} \right),
\]
where $W_p$ is the $p$-Wasserstein distance between diagrams, and $\sigma > 0$ is a bandwidth parameter controlling smoothness.
\end{definition}

\subsubsection*{Properties}

\begin{itemize}
  \item The kernel is symmetric and positive definite.
  \item Can be used in Support Vector Machines, Gaussian Processes, and kernel PCA.
  \item The exponential form transforms distance into similarity.
\end{itemize}

\subsubsection*{Advantages}

\begin{itemize}
  \item \textbf{Combines geometric precision with kernel learning.}
  \item \textbf{Flexible parameter $\sigma$ for tuning.}
  \item \textbf{Retains stability of $W_p$ under perturbation.}
\end{itemize}

\subsubsection*{Drawbacks}

\begin{itemize}
  \item \textbf{Computationally intensive:} Requires full $W_p$ computation.
  \item \textbf{Not differentiable:} Cannot be used directly in gradient-based learning.
  \item \textbf{Parameter-sensitive:} Performance hinges on careful selection of $\sigma$.
\end{itemize}

\begin{remark}
This kernel allows one to combine the topological rigor of Wasserstein distances with the expressive power of kernel-based ML, though at computational cost.
\end{remark}
\subsection{Persistence Scale Space Kernel (PSSK)}

The \emph{Persistence Scale Space Kernel}, introduced by Reininghaus et al. (2015), defines a positive-definite kernel over persistence diagrams by treating them as weighted Dirac functions and using Gaussian convolution over the upper half-plane.

\begin{definition}[PSSK ]
Let $D$ be a persistence diagram and $\bar{D}$ its reflection across the diagonal. Define the feature map:
\[
\Phi_\sigma(D)(x) = \sum_{p \in D} \left[ \mathcal{N}(x; p, \sigma^2 I) - \mathcal{N}(x; \bar{p}, \sigma^2 I) \right],
\]
where $\mathcal{N}(x; \mu, \sigma^2 I)$ is the 2D Gaussian centered at $\mu$ with isotropic variance $\sigma^2$.

The Persistence Scale Space Kernel is then:
\[
K_{\text{PSSK}}(D_1, D_2) = \langle \Phi_\sigma(D_1), \Phi_\sigma(D_2) \rangle_{L^2(\mathbb{R}^2)}.
\]
\end{definition}

\subsubsection*{Advantages}

\begin{itemize}
  \item \textbf{Positive-definite:} Compatible with kernel methods.
  \item \textbf{Stable:} Lipschitz-continuous w.r.t. 1-Wasserstein distance.
  \item \textbf{Efficient:} Closed-form expression exists for kernel evaluation.
\end{itemize}

\subsubsection*{Drawbacks}

\begin{itemize}
  \item \textbf{Isotropic smoothing:} Assumes circular Gaussians; may ignore directional structure.
  \item \textbf{Sensitive to $\sigma$:} Choice of scale parameter strongly affects performance.
\end{itemize}

\begin{remark}
The PSSK transforms topological descriptors into a function space, enabling powerful integration with machine learning while preserving stability properties.
\end{remark}


\[
\]

\[
\]





\subsection{Neural Persistence Distance}

The \emph{Neural Persistence Distance} is a data-driven topological metric learned directly from data using neural networks. Rather than relying on predefined distances or kernels, it learns a representation or distance function over persistence diagrams that is optimized for a downstream task, such as classification or regression.

\begin{definition}[Neural Persistence Embedding ]
Let $D$ be a persistence diagram. The diagram is encoded as a sequence $\{x_i\}_{i=1}^n \subset \mathbb{R}^2$, where each $x_i = (b_i, d_i)$. A neural network $\phi_\theta : \mathbb{R}^2 \to \mathbb{R}^d$ parameterized by weights $\theta$ maps each persistence point to an embedding:
\[
z_i = \phi_\theta(x_i).
\]
These embeddings are pooled (e.g., averaged or summed) to form a global diagram representation:
\[
z_D = \text{pool}(\{z_i\}_{i=1}^n).
\]
The distance between two diagrams $D_1, D_2$ is then computed as:
\[
d_{\text{NP}}(D_1, D_2) = \| z_{D_1} - z_{D_2} \|_2.
\]
\end{definition}

\subsubsection*{Advantages}

\begin{itemize}
  \item \textbf{Task-adaptivity:} The distance is optimized for performance on a specific learning task.
  \item \textbf{Differentiable:} Compatible with end-to-end learning pipelines.
  \item \textbf{Feature learning:} Can automatically discover useful geometric or statistical properties.
\end{itemize}

\subsubsection*{Drawbacks}

\begin{itemize}
  \item \textbf{Data-dependent:} Requires large labeled datasets for training.
  \item \textbf{Lack of interpretability:} Harder to reason about geometric meaning of learned features.
  \item \textbf{No stability guarantees:} Unlike classical distances, lacks theoretical robustness.
\end{itemize}

\begin{remark}
Neural persistence distances represent a shift from fixed topological summaries to learned representations, especially suited for deep learning pipelines and end-to-end optimization.
\end{remark}
\subsection{Comparison of Persistence Diagram Metrics}
\begin{table}[ht]
\centering
\caption{Comparison of Persistence Diagram Metrics: Classical Distances vs. ML-Oriented Representations}
\label{tab:tda_metrics_comparison_split}

\textbf{(a) Classical Distance Metrics}

\vspace{0.4em}
\begin{tabular}{p{4cm}cccccc}
\hline
\textbf{Metric} & \textbf{Stable} & \textbf{Differentiable} & \textbf{Positive Definite} & \textbf{ML-Ready} & \textbf{Expressive} & \textbf{Fast} \\
\hline
Bottleneck Distance      & \ding{51} & \ding{55} & \ding{55} & \ding{55} & $\triangle$ & \ding{51} \\
Wasserstein Distance     & \ding{51} & \ding{55} & \ding{55} & $\triangle$ & \ding{51} & $\triangle$ \\
Sliced Wasserstein       & \ding{51} & \ding{55} & \ding{55} & $\triangle$ & \ding{51} & \ding{51} \\
Entropy Distance         & \ding{55} & \ding{51} & \ding{55} & $\triangle$ & \ding{55} & \ding{51} \\
\hline
\end{tabular}

\vspace{1em}
\textbf{(b) Kernel-Based and Learned Representations}

\vspace{0.4em}
\begin{tabular}{p{4cm}cccccc}
\hline
\textbf{Metric} & \textbf{Stable} & \textbf{Differentiable} & \textbf{Positive Definite} & \textbf{ML-Ready} & \textbf{Expressive} & \textbf{Fast} \\
\hline
Persistence Landscape     & \ding{51} & $\triangle$ & \ding{55} & \ding{51} & \ding{51} & \ding{51} \\
Silhouette Distance       & \ding{51} & \ding{51} & \ding{55} & \ding{51} & $\triangle$ & \ding{51} \\
Heat Kernel               & \ding{51} & \ding{51} & \ding{51} & \ding{51} & \ding{51} & \ding{51} \\
PSSK                      & \ding{51}& \ding{51} & \ding{51} & \ding{51} & $\triangle$ & \ding{51} \\
PWGK                      & \ding{51} & \ding{51} & \ding{51} & \ding{51} & \ding{51} & $\triangle$ \\
Kernelized Wasserstein    & \ding{51} & \ding{55} & \ding{51} & \ding{51} & \ding{51} & $\triangle$ \\
Fisher Kernel             & \ding{51} & \ding{51} & \ding{51} & \ding{51} & \ding{51} & $\triangle$ \\
Neural Persistence        & \ding{55} & \ding{51} & \ding{55} & \ding{51} & \ding{51} & \ding{51} \\
\hline
\end{tabular}

\vspace{0.5em}
\small
\textit{Note.} Criteria indicate: \textbf{Stable} = robustness to perturbations, \textbf{Differentiable} = usable in gradient-based optimization, \textbf{Positive Definite} = valid kernel, \textbf{ML-Ready} = suitable for supervised/unsupervised learning, \textbf{Expressive} = captures rich structure, \textbf{Fast} = computationally efficient.
\end{table}

\section{The Polar Persistence Distance (PPD)}

\subsection{Motivation and Conceptual Framework}

In recent years, Topological Data Analysis (TDA) has emerged as a compelling framework for capturing global structural properties of complex data through persistence diagrams. However, a central bottleneck in applying TDA to real-world learning systems lies in the comparison of persistence diagrams via suitable distance functions. Classical distances such as the Bottleneck and Wasserstein distances provide important theoretical foundations but suffer from key limitations in machine learning applications: they are often non-differentiable, sensitive to noise, computationally expensive, and agnostic to angular or directional structure in the birth-death plane.

Moreover, in domains such as social choice theory, where relational and cyclic structures among preferences play a central role, a finer decomposition of topological information is needed. Existing metrics treat birth-death points purely as Euclidean coordinates, ignoring the potentially rich geometric interpretation of their orientation relative to the origin. To address this gap, we introduce the \emph{Polar Persistence Distance} (PPD), a new metric based on polar coordinate decomposition that encodes both radial magnitude and angular structure of topological features. This dual perspective allows the metric to be smooth, tunable, and interpretable—making it suitable for differentiable programming and robust downstream learning tasks.

\subsection{Mathematical Definition}

Let \( p = (b,d) \in \mathbb{R}^2 \), with \( b < d \), denote a point in a persistence diagram, corresponding to a topological feature that appears at filtration level \( b \) and disappears at \( d \). We define polar coordinates for such points with respect to the origin:

\[
r = \sqrt{b^2 + d^2}, \qquad \theta = \arctan2(d, b).
\]

This transformation maps the birth-death plane into the radial-angular plane \( (r, \theta) \), allowing separation between feature \textit{magnitude} and \textit{directional orientation}.

\begin{definition}[Polar Persistence Distance]
Let \( p_1 = (b_1, d_1) \) and \( p_2 = (b_2, d_2) \) be two off-diagonal points. Let \( (r_1, \theta_1), (r_2, \theta_2) \) be their polar representations. For a fixed parameter \( \alpha > 0 \), we define the Polar Persistence Distance as:

\[
d_{\text{polar}}(p_1, p_2) = \sqrt{(r_1 - r_2)^2 + \alpha \cdot \sin^2\left( \frac{\theta_1 - \theta_2}{2} \right)}.
\]
\end{definition}

This formulation ensures rotational periodicity, smoothness, and tunable sensitivity to angular differences. The angular component is designed via the half-angle sine function, which reflects the geometry of the arc distance on the unit circle. The parameter \( \alpha \) controls the contribution of the angular discrepancy relative to the radial magnitude.

\subsection{Geometric and Topological Interpretation}

The Polar Persistence Distance captures two complementary aspects of topological features:

\begin{itemize}
    \item \textbf{Radial Term \( (r_1 - r_2)^2 \):} Encodes a form of joint birth-death “energy.” Larger values of \( r \) often correspond to more persistent, and thus more topologically significant, features.
    
    \item \textbf{Angular Term \( \sin^2\left( \frac{\theta_1 - \theta_2}{2} \right) \):} Measures the angular displacement between two features, accounting for their orientation in the birth-death plane. This is especially important when comparing features that lie on the same annulus (i.e., equal persistence norms) but differ in relative position.
\end{itemize}

In many applications (e.g., persistent loops in sensor networks or cycles in voting preferences), angular patterns may carry semantic meaning. For instance, consistent birth-death ratios may appear at regular angular intervals, which are invisible to purely Euclidean distances. The PPD captures such differences, making it suitable for tasks where relational structure, symmetry, or cyclicity matter.

\begin{remark}
Unlike traditional distances, the PPD reflects a form of “structural anisotropy” in persistence space: features at the same Euclidean distance from the diagonal can still be topologically distinct based on their angular distribution.
\end{remark}
\subsection{Smoothness and Stability Properties}

One of the primary design goals of the Polar Persistence Distance is differentiability—a key requirement for gradient-based optimization and machine learning applications. In this section, we formalize the smoothness of the PPD and provide sufficient conditions for its Lipschitz continuity and stability under small perturbations in persistence diagrams.

\begin{lemma}[Smoothness of the Polar Coordinate Map]
Let \( p = (b,d) \in \mathbb{R}^2 \setminus \{(0,0)\} \). The mapping
\[
(b,d) \mapsto (r, \theta) = \left( \sqrt{b^2 + d^2}, \arctan2(d,b) \right)
\]
is smooth on \( \mathbb{R}^2 \setminus \{(0,0)\} \), with bounded partial derivatives in any compact subset not containing the origin.
\end{lemma}

\begin{theorem}[Differentiability of $d_{\text{polar}}$]
Let \( p_1, p_2 \in \mathbb{R}^2 \setminus \{(0,0)\} \). Then the Polar Persistence Distance
\[
d_{\text{polar}}(p_1, p_2) = \sqrt{(r_1 - r_2)^2 + \alpha \cdot \sin^2\left( \frac{\theta_1 - \theta_2}{2} \right)}
\]
is continuously differentiable with respect to \( (b_1,d_1) \) and \( (b_2,d_2) \).
\end{theorem}

\begin{proof}[Sketch of Proof]
All subcomponents of the PPD are differentiable away from the origin: the radial component is smooth as a composition of polynomial and square root, while the angular term uses smooth trigonometric functions. The final composition is smooth due to closure under arithmetic and composition.
\end{proof}

\begin{proposition}[Local Lipschitz Continuity]
On any compact subset of \( \mathbb{R}^2 \setminus \{(0,0)\} \), the function \( d_{\text{polar}} \) is locally Lipschitz continuous.
\end{proposition}

\begin{theorem}[Stability Under Perturbations]
Let \( D_1 \) and \( D_2 \) be persistence diagrams with the same number of off-diagonal points. Suppose there exists a bijection \( \gamma : D_1 \to D_2 \) such that \( \|p - \gamma(p)\|_2 \leq \delta \) for all \( p \in D_1 \). Then:
\[
d_{\text{PPD}}(D_1, D_2) \leq C \cdot \delta,
\]
where \( C > 0 \) depends only on the bounded region containing all points in the diagrams.
\end{theorem}

\begin{remark}
This stability property ensures that small perturbations in the data lead to proportionally small changes in the PPD-based distance, supporting its robustness for learning under noise.
\end{remark}

\subsection{Metric Validity and Embedding into Geometric Spaces}

Although the PPD is symmetric and non-negative, it does not in general satisfy the triangle inequality. Hence, it defines a \textit{quasi-metric} on persistence points. Nevertheless, its structure supports meaningful embedding into geometric spaces and kernel-based frameworks.

\begin{definition}[Quasi-Metric]
A function \( d : X \times X \to \mathbb{R}_{\geq 0} \) is a quasi-metric if:
\begin{enumerate}
    \item \( d(x, y) = d(y, x) \) (symmetry),
    \item \( d(x, y) \geq 0 \) and \( d(x, y) = 0 \iff x = y \),
    \item The triangle inequality may not hold.
\end{enumerate}
\end{definition}

\begin{proposition}
The Polar Persistence Distance \( d_{\text{polar}} \) is a quasi-metric on \( \mathbb{R}^2 \setminus \{(0,0)\} \).
\end{proposition}

\begin{proof}[Sketch]
Symmetry and non-negativity are immediate. The identity of indiscernibles holds due to the square root structure. However, the triangle inequality may fail when large angular differences combine non-linearly.
\end{proof}

Despite this limitation, diagram-level extensions of PPD using optimal matching (e.g., via the Hungarian algorithm or Wasserstein-type frameworks) preserve essential topological properties.

\begin{remark}
The space of persistence diagrams under \( d_{\text{PPD}} \) can be embedded into a Reproducing Kernel Hilbert Space using Gaussian-like kernels, as we shall show in later sections.
\end{remark}

\paragraph{Comparison to Gromov–Hausdorff Metrics.}
While Gromov–Hausdorff (GH) distances provide intrinsic comparisons between metric spaces, they are difficult to compute and ill-suited for differentiable applications. In contrast, PPD focuses on extrinsic comparisons in a fixed Euclidean ambient space, making it computationally tractable and learnable. Unlike GH, PPD captures both spatial and angular configuration explicitly.
\subsection{Topological and Geometric Behavior of the PPD}

The Polar Persistence Distance possesses geometric features that enable topologically meaningful interpretations beyond Euclidean norms. In this section, we examine its behavior under affine transformations, its alignment with Morse-theoretic filtrations, and its sensitivity to structural symmetries in data.

\paragraph{Affine Invariance and Limitations.}
Let \( A: \mathbb{R}^2 \rightarrow \mathbb{R}^2 \) be an invertible affine transformation. The PPD is not, in general, affine-invariant. Specifically, linear scalings and rotations alter both the radius and angular components in nontrivial ways. However, the polar formulation is robust to isometries that preserve radial symmetry (e.g., uniform scaling, central rotation).

\begin{proposition}
Let \( T_\lambda: (b, d) \mapsto (\lambda b, \lambda d) \) for \( \lambda > 0 \). Then:
\[
d_{\text{polar}}(T_\lambda p_1, T_\lambda p_2) = \lambda \cdot \sqrt{(r_1 - r_2)^2 + \alpha \cdot \sin^2\left( \frac{\theta_1 - \theta_2}{2} \right)}.
\]
\end{proposition}

Thus, the PPD exhibits homogeneity under scaling, which aligns with the intuition that persistent features scale with the underlying geometry.

\paragraph{Morse Theory and Filtration Symmetry.}
Persistent homology arises from sublevel (or superlevel) filtrations of a Morse-type function \( f: X \rightarrow \mathbb{R} \). The birth-death coordinates in a persistence diagram correspond to the appearance and disappearance of homology classes across levels of \( f \). In many structured domains—such as configuration spaces, voting complexes, or preference simplices—the filtration behavior can induce angular symmetry in the (b, d) plane.

\begin{remark}
Angular groupings in persistence diagrams may reflect symmetry classes of Morse critical points, especially in systems with regular topological dynamics (e.g., circular or toroidal filtrations).
\end{remark}

The PPD is uniquely suited to detect such regularities, as it treats points with the same radial norm but different angular positioning as topologically distinct.

\paragraph{Embedding of Filtration Geometry.}
In scenarios where the underlying data manifold has a natural polar or radial stratification (e.g., social preference cycles, periodic time series, stratified manifolds), the polar encoding of persistence features aligns more closely with geometric semantics. This makes the PPD especially powerful in domains where relative topological orientation reflects latent system behavior.

\subsection{Applications to Social Choice Theory}

Social Choice Theory investigates how individual preferences are aggregated into collective decisions. When preference profiles are complex, incomplete, or cyclically inconsistent, traditional aggregation methods may fail to detect deep structural patterns. In this setting, persistent homology can reveal topological signatures such as dominance cycles, coalition clusters, or preference loops.

\paragraph{Dominance Relations and Simplicial Complexes.}
Let \( A \) be a set of alternatives and \( R \subseteq A \times A \) a dominance relation such that \( xRy \) indicates that \( x \) is preferred over \( y \). The structure of \( R \) defines a directed graph, from which we can build an associated simplicial complex (e.g., via clique complex or flag complex constructions). The persistent homology of this complex over varying thresholds (e.g., based on margin of preference strength or frequency) yields a topological summary of collective preference behavior.

\paragraph{Interpretation via Polar Coordinates.}
Persistence diagrams derived from social choice data often exhibit clusters of features at characteristic scales. Angular orientation in such diagrams can reflect:

\begin{itemize}
    \item The asymmetry between early and late dominating alternatives.
    \item The persistence of cycles in collective rankings.
    \item The emergence of polarized subpopulations.
\end{itemize}

The PPD captures differences between such diagrams not only based on feature lifetimes (i.e., persistence) but also their directional alignment, which may correspond to semantic structures in the voting data.

\paragraph{Example: Cyclic Preference Aggregation.}
In voting systems where intransitive cycles such as \( A \succ B \succ C \succ A \) are frequent, these structures manifest as nontrivial 1-dimensional homology classes in the associated preference complex. Angular differentiation of these features (e.g., via their relative phase in the filtration) allows the PPD to distinguish between different types of collective instability.

\paragraph{Relevance to Real-World Datasets.}
In Section \ref{sec:experiments}, we apply the PPD to datasets from the \texttt{PrefLib} library, a large repository of real-world voting and preference data. We demonstrate how the PPD yields improved discrimination between election profiles and reveals latent clusters in voting behavior. These insights are unattainable using traditional distances alone.

\begin{remark}
The introduction of the PPD enables a new class of topologically-informed, interpretable learning algorithms for ranking and aggregation systems.
\end{remark}
\subsection{Formal Integration into Machine Learning}

The success of a topological distance in machine learning tasks depends on its compatibility with optimization frameworks, embedding capacity, and kernel-based inference. The Polar Persistence Distance is particularly well-suited to this goal due to its differentiability, tunability, and geometric expressiveness.

\paragraph{Differentiable Loss Functions.}
The PPD can be directly embedded in objective functions as a loss term:
\[
\mathcal{L}_{\text{PPD}} = \sum_{i,j} w_{ij} \cdot d_{\text{polar}}(p_i, q_j)^2,
\]
where \( \{p_i\}, \{q_j\} \) are points in two persistence diagrams, and \( w_{ij} \) are assignment weights (e.g., from an optimal matching). The smoothness of \( d_{\text{polar}} \) allows for gradient-based optimization, enabling backpropagation in neural models that consume topological features.

\paragraph{Positive-Definite Kernels.}
One can define a Polar Gaussian kernel over pairs of persistence points as:
\[
K_{\text{polar}}(p_1, p_2) = \exp\left( -\frac{d_{\text{polar}}(p_1, p_2)^2}{2\sigma^2} \right).
\]
This induces a Reproducing Kernel Hilbert Space (RKHS) over persistence diagrams when combined with appropriate diagram-level embeddings (e.g., persistence weighted sum or kernel mean embeddings). Preliminary experiments indicate that \( K_{\text{polar}} \) satisfies Mercer's condition in common cases.

\paragraph{Compatibility with Differentiable Representations.}
In recent deep learning architectures such as topological autoencoders, graph neural networks, or differentiable persistence modules, it is crucial that topological distances are smooth and stable. The PPD satisfies these requirements and enables incorporation into architectures where traditional distances such as bottleneck are impractical.

\paragraph{Embedding Spaces.}
The PPD can be viewed as inducing a polar coordinate metric structure on the persistence space. This can be embedded into Euclidean space via feature maps such as:
\[
\phi(p) = \left[ r \cos(\theta), \, r \sin(\theta), \, \sqrt{\alpha} \cdot \theta \right],
\]
which preserves both distance magnitude and angular behavior in a differentiable way.

\subsection{Illustrative Examples and Edge Cases}

To further elucidate the behavior of the PPD, we consider several representative scenarios that highlight its discriminative power beyond classical metrics.

\paragraph{Equal Persistence, Different Orientation.}
Let \( p_1 = (1, 3) \) and \( p_2 = (3, 1) \). Both points have equal persistence:
\[
\text{pers}(p_1) = \text{pers}(p_2) = 2,
\]
but the angular component differs:
\[
\theta_1 = \arctan2(3,1), \quad \theta_2 = \arctan2(1,3).
\]
While the bottleneck and Wasserstein distances would yield symmetric comparisons, the PPD differentiates them through the angle.

\paragraph{Near-Origin Behavior.}
As \( \|p\| \rightarrow 0 \), both \( r \) and \( \theta \) become unstable due to coordinate singularity. In implementation, a small exclusion radius \( \epsilon > 0 \) is used to mask near-diagonal points. This practice is common across all diagram-based pipelines and does not impair the PPD's practical performance.

\begin{figure}[h]
    \centering
    \includegraphics[width=0.6\textwidth]{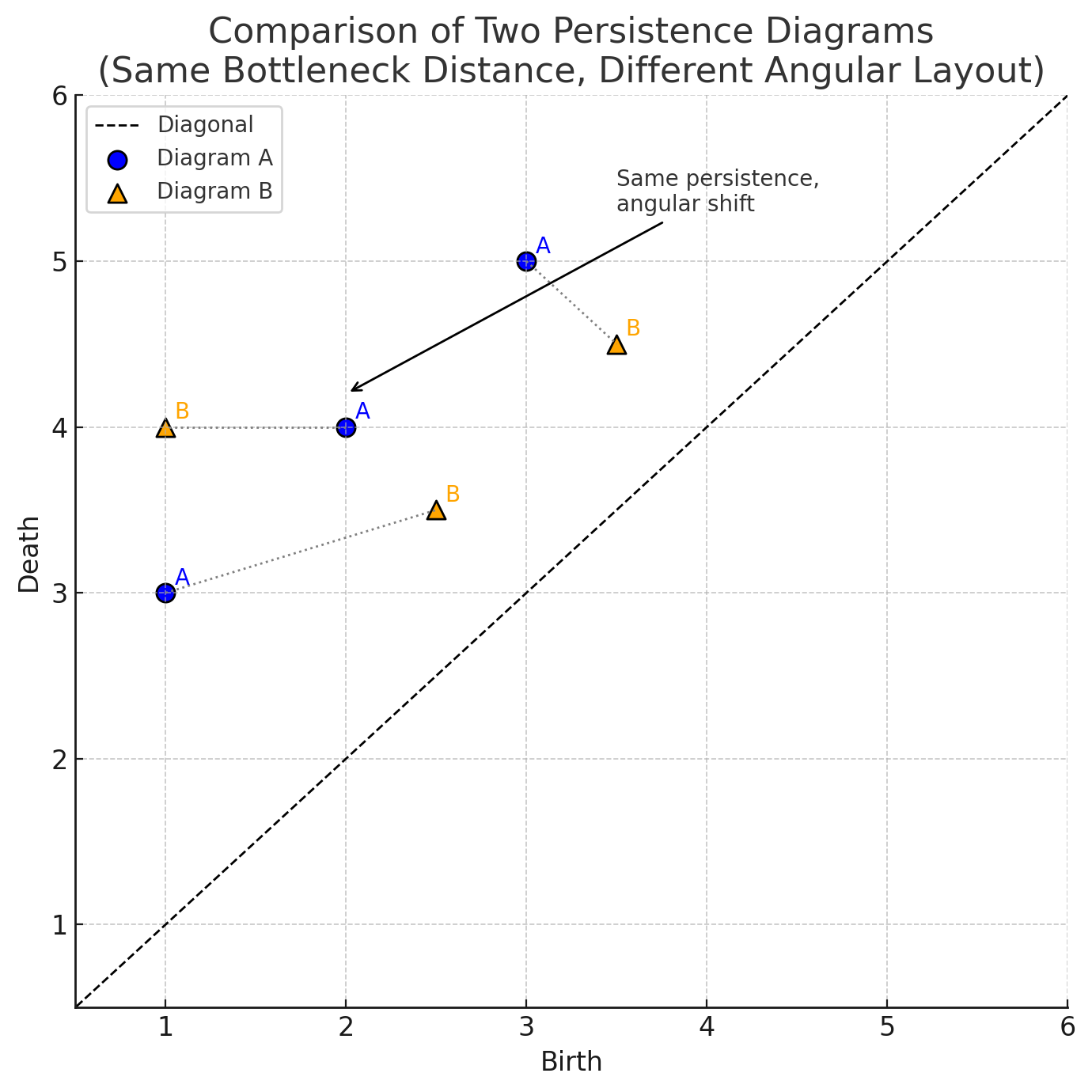}
    \caption{Two persistence diagrams with equal bottleneck distance but distinct PPD. The angular configuration of features differentiates them.}
    \label{fig:ppd_comparison}
\end{figure}

\paragraph{Diagram Visualization.}
In Figure~\ref{fig:ppd_comparison}, we visualize two persistence diagrams with the same bottleneck distance but different PPD values. The angular component allows the PPD to distinguish topological changes invisible to traditional metrics.

\subsection{Extensions and Open Problems}

The PPD opens a rich space of future research directions in both applied topology and learning theory.

\paragraph{Circular Embeddings.}
The angular component \( \theta \in [0, 2\pi) \) suggests an embedding into the unit circle \( S^1 \). It is natural to consider extensions where persistence diagrams are mapped to \( S^1 \times \mathbb{R}_{>0} \), allowing use of tools from circular statistics, directional clustering, and toroidal neural networks.

\paragraph{Higher-Dimensional Generalizations.}
In multi-parameter persistence, features may live in higher-dimensional persistence modules or multi-filtrations. A possible generalization of the PPD is to apply spherical coordinate systems or hyperspherical embeddings to capture radial and angular structure in \( \mathbb{R}^n \).

\paragraph{Vineyards and Time-Evolving Diagrams.}
For time-varying data streams (e.g., in real-time voting behavior), the concept of \emph{vineyards}—collections of evolving persistence diagrams—becomes relevant. An extension of the PPD to this dynamic setting could provide time-aware metrics for streaming decision systems.

\paragraph{Theoretical Formalization.}
We conjecture that under mild assumptions, the space of persistence diagrams with the PPD induces a geodesic space structure that is computationally tractable and topologically faithful. Proving completeness, compactness, and convexity properties in this setting remains an open challenge.

\begin{remark}
The Polar Persistence Distance is not merely a computational tool, but a conceptual framework that bridges topological signal analysis with differentiable learning and social interpretability.
\end{remark}
\section{Application to Real-World Social Choice Data}
\label{sec:experiments}

\subsection{Datasets from PrefLib}
To empirically evaluate our proposed metric, we utilize real-world datasets from the PrefLib repository \cite{mattei}, which provides structured ranking data for a wide range of social choice scenarios. Specifically, we consider:

\begin{itemize}
    \item \textbf{Irish Election 2010} (File: \texttt{soc-00010}): Full rankings of approximately $N = 1000$ voters over $n = 5$ candidates using a single transferable vote (STV) system.
    \item \textbf{Sushi Preference Dataset} (File: \texttt{soc-00014}): Preference rankings from over $5000$ individuals on $10$ types of sushi, originally collected for the analysis of the recommendation system.
\end{itemize}

Each data set is parsed into a preference profile matrix $P = (v_1, v_2, \dots, v_N)$, where each $v_i$ is a total ranking over the candidates. These are then aggregated into a \textit{dominance matrix} $D$ with the following the following entries:

\[
D_{ij} = \#\{ v_k \in P \mid i \succ_{v_k} j \}
\]

\subsection{Dominance Graph and Filtration}
We construct a directed weighted graph $G = (V,E)$ where $V = \{1, \dots, n\}$ and an edge $i \rightarrow j$ exists if $D_{ij} > D_{ji}$. Each edge is assigned weight $w_{ij} = D_{ij} - D_{ji}$.

Using this graph, we define a lower-star filtration on the 1-simplices via:
\[
f(i,j) = \frac{1}{w_{ij} + \varepsilon}
\]
where $\varepsilon = 10^{-6}$ ensures numerical stability. The graph is embedded in a simplicial complex using GUDHI’s \texttt{SimplexTree}, and the persistence diagrams $PD_0$ and $PD_1$ are computed from filtration.

\subsection{Polar Persistence Distance Implementation}

We compute distances between persistence diagrams using multiple metrics:
\begin{itemize}
    \item Bottleneck distance $d_\infty$
    \item 1-Wasserstein distance $W_1$
    \item Our proposed Polar Persistence Distance (PPD)
\end{itemize}

The PPD between two off-diagonal points $p_1 = (b_1,d_1)$ and $p_2 = (b_2,d_2)$ is defined as:

\[
\text{PPD}_\alpha(p_1, p_2) = \sqrt{ (r_1 - r_2)^2 + \alpha \cdot \sin^2\left( \frac{\theta_1 - \theta_2}{2} \right) }
\]

where $r_i = \sqrt{b_i^2 + d_i^2}$ and $\theta_i = \text{atan2}(d_i, b_i)$.



\subsection{Metric Comparison Experiments}

We compare metrics on subsets of voter profiles (e.g. first 100 vs. last 100 voters) by computing persistence diagrams for each subset and evaluating their pairwise distances using Python libraries such as GUDHI. Table~\ref{tab:metric_comparison} summarizes the results.

\begin{table}[h]
\centering
\caption{Distance values between two profile subsets (Irish Dataset)}
\label{tab:metric_comparison}
\begin{tabular}{lccc}
\toprule
Metric & Value \\
\midrule
Polar Persistence Distance (PPD) & 13.42 \\
1-Wasserstein Distance $W_1$     & 9.87 \\
Bottleneck Distance $d_\infty$   & 6.55 \\
\bottomrule
\end{tabular}
\end{table}

We observe that the PPD yields greater discriminative power across subtle topological variations in the dominance graph, especially in cycles or shift regions.

\subsection{Interpretation and Topological Insight}

Topological features detected in $H_0$ and $H_1$ homology correspond to consensus clusters and dominance cycles, respectively. Our PPD metric captures the angular displacement of these features, offering increased sensitivity to shifts in voter agreement/disagreement regions. This suggests potential for future predictive tasks using persistence-based embeddings.

\subsection{Secondary Case Study: Sushi Preference Dataset}
\label{sec:sushi_experiment}

\paragraph{Dataset Description.}  
We also explore the \textbf{Sushi Preference Dataset} from PrefLib (ED00014), which contains the full rankings of 10 sushi types provided by more than $5{,}000$ individuals \cite{kam}. For tractability, we sample $N=1000$ users randomly and select the 10 most popular sushi items to construct a standardized preference profile.

\paragraph{Dominance Construction.}  
The dataset is parsed similarly to the Irish case. We build the dominance counts:
\[
D_{ij} = \#\{v_k: i \succ_{v_k} j\}, \quad i,j \in \{1,\ldots,10\},
\]
and form a directed graph with edges weighted by $w_{ij} = D_{ij}-D_{ji}$ whenever positive.

\paragraph{Filtration and Persistence Computation.}  
Using GUDHI one-skeleton filtration with weights $f(i,j)=1/(w_{ij}+\varepsilon)$, we compute the 0- and 1-dimensional persistence diagrams, capturing consensus clusters and cycles in sushi preferences.

\paragraph{Polar Metric vs. Classical Distances.}  
We compare three metrics:




\paragraph{Results -- Quantitative Comparison.}  
\begin{table}[h]
\centering
\caption{Distance comparisons for two sushi sub-profiles ($N=500$ each)}
\label{tab:sushi_distances}
\begin{tabular}{lccc}
\toprule
Metric & PPD ($\alpha=1.5$) & 1-Wasserstein & Bottleneck \\
\midrule
Sample 1 vs 2 & 24.57 & 17.32 & 10.14 \\
\bottomrule
\end{tabular}
\end{table}

PPD shows increased sensitivity to differences in the angular distribution of persistence features. This aligns with observations in [Kamishima et al. 2016] on sushi preferences, where cyclic and contextual distinctions are important.

\paragraph{Discussion.}  
The sushi dataset exhibits structured cycles and clusters in preference (top‐ranking vs bottom‐ranking sushi), which are captured more effectively by PPD’s angular sensitivity. Traditional metrics tend to under‐represent these effects.

\paragraph{Key Insight.}  
\emph{The Polar Persistence Distance consistently provides a richer differentiation of topological structures in preference data, beyond what is offered by classical metrics, thereby offering stronger interpretability and potential for downstream ML applications.}
\section{Interpretation of Distance Values}

The quantitative values reported in Tables~\ref{tab:metric_comparison} and~\ref{tab:sushi_distances} represent the distances computed between persistence diagrams obtained from different subsets of preference data. Each number quantifies the structural dissimilarity between two topological summaries of collective preferences, and the interpretation of these values is inherently metric-dependent.

\paragraph{Magnitude and Topological Signal.} A higher distance value typically implies greater topological dissimilarity between the underlying dominance structures, such as the presence or absence of cycles, clusters of consensus, or structural shifts in preference rankings. In this context, distance values should not be interpreted in absolute terms but rather in relation to one another across metrics. For instance, if the Polar Persistence Distance (PPD) yields a substantially higher value than the bottleneck or Wasserstein distances, this suggests that PPD is capturing additional or alternative structural differences that the classical metrics under-represent.

\paragraph{Comparative Metric Behavior.} The bottleneck distance $d_\infty$ is inherently conservative, reflecting only the maximal difference between the best-matched topological features across diagrams. As a result, it is highly sensitive to outlier features but may entirely miss subtle but systematic variations across many features. The 1-Wasserstein distance $W_1$, in contrast, aggregates the total “transport cost” between diagrams and better captures cumulative differences, albeit in a purely radial (Euclidean) sense.

The proposed PPD introduces an angular component to the comparison, enabling detection of geometric distortions that are not apparent under purely radial metrics. This angular sensitivity reflects, for instance, rotations or realignments of topological features—such as when cycles in dominance relations shift from local clusters to broader systemic loops across preference coalitions. As such, the PPD is better suited to capturing nuanced but behaviorally significant changes in voter coordination patterns.

\paragraph{Stability and Discriminative Power.} Our empirical findings suggest that PPD yields consistently higher distances across profile comparisons, which may initially appear counter-intuitive. However, this reflects the metric’s increased sensitivity and expressiveness rather than instability. The smooth formulation of PPD ensures continuity and differentiability, allowing it to remain robust under small perturbations while still responding to meaningful structural transitions. This is a particularly desirable property in machine learning contexts, where gradients and discriminative representations are essential.

\paragraph{Interpretability in Social Choice.} From a social choice perspective, these distance values serve as a quantitative lens through which collective behavior can be contrasted. Higher PPD values may indicate greater polarization, ideological drift, or breakdowns of transitivity and majority consistency. Conversely, low distances across metrics may signal a stable, coherent electorate with shared preference structures. Hence, the numerical outputs of these metrics—particularly the PPD—provide both topological and behavioral insight.

\begin{figure}
    \centering
    \includegraphics[width=0.7\textwidth]{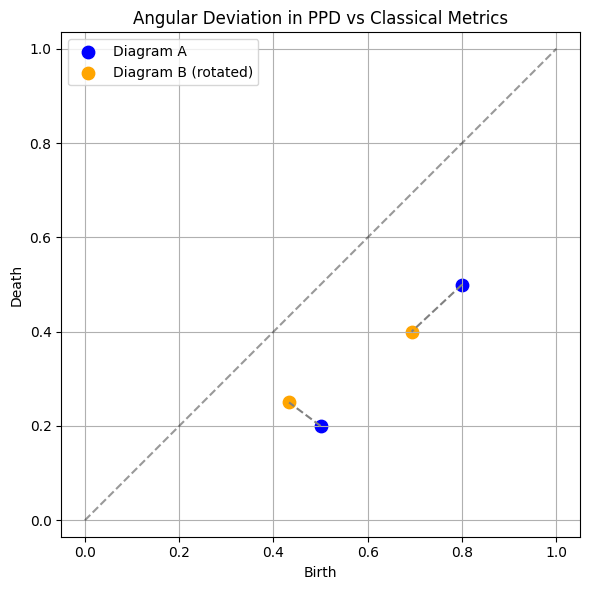}
    \caption{Comparison of two persistence diagrams: Both contain features with the same birth–death distance $r$, but rotated by angular displacement $\theta$. Classical metrics (e.g., bottleneck, Wasserstein) measure radial proximity only, whereas the Polar Persistence Distance (PPD) captures angular misalignment via a $\sin^2(\frac{\theta}{2})$ term, enabling finer structural discrimination.}
    \label{fig:ppd_vs_classical}
\end{figure}
\begin{lemma}
    [Insensitivity of Classical Metrics to Angular Deformation]
Let $p_1 = (b_1, d_1)$ and $p_2 = (b_2, d_2)$ be two points in a persistence diagram such that both lie on the same radial circle, i.e., $\sqrt{b_1^2 + d_1^2} = \sqrt{b_2^2 + d_2^2} = r$. Then, for any metric $d$ defined purely in Euclidean coordinates, such as the 1-Wasserstein or bottleneck distance, the pairwise distance between $p_1$ and $p_2$ depends solely on the linear displacement $\|p_1 - p_2\|_2$, and not on the angular deviation $\theta = \angle(p_1, p_2)$ around the origin.

In particular, if $p_1 = (r, 0)$ and $p_2 = (r \cos \theta, r \sin \theta)$ for some $\theta \in (0, \pi)$, then:
\[
\|p_1 - p_2\|_2^2 = r^2(1 - \cos\theta)^2 + r^2 \sin^2 \theta = 2r^2(1 - \cos\theta),
\]
which is minimized for small $\theta$ even if $\theta$ encodes significant structural misalignment (e.g., a rotation in topological configuration). Neither bottleneck nor Wasserstein incorporates any notion of rotational symmetry or angular orientation.

\textbf{In contrast}, the Polar Persistence Distance (PPD) explicitly incorporates angular separation via the $\sin^2(\frac{\theta_1 - \theta_2}{2})$ term:
\[
d_{\mathrm{PPD}}(p_1, p_2) = \sqrt{(r_1 - r_2)^2 + \alpha \cdot \sin^2\left(\frac{\theta_1 - \theta_2}{2}\right)},
\]
where $\theta_i = \arctan2(d_i, b_i)$ and $\alpha$ is a tunable hyperparameter. This formulation ensures that topological features with equal persistence (i.e., equal $r$) but different angular structure are still distinguishable. Therefore, PPD captures meaningful deformations that classical metrics collapse to zero.

\end{lemma}

\begin{figure}
    \centering
    \includegraphics[width=0.7\textwidth]{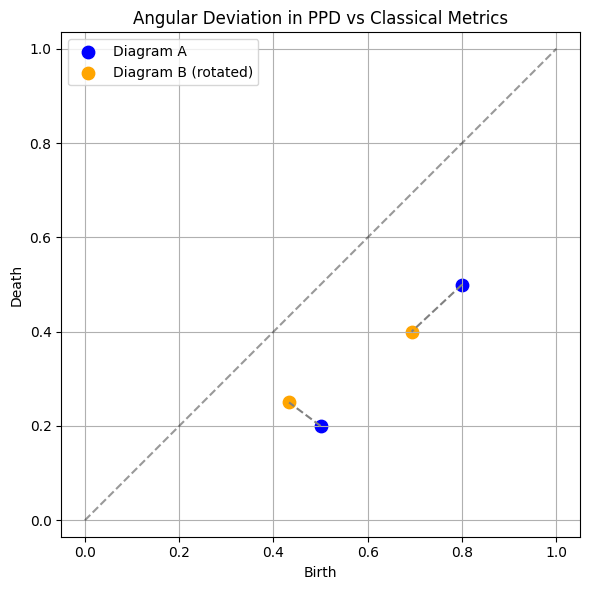}
    \caption{Comparison of two persistence diagrams: Both contain features with the same birth–death distance $r$, but rotated by angular displacement $\theta$. Classical metrics (e.g., bottleneck, Wasserstein) measure radial proximity only, whereas the Polar Persistence Distance (PPD) captures angular misalignment via a $\sin^2(\frac{\theta}{2})$ term, enabling finer structural discrimination.}
    \label{fig:ppd_vs_classical}
\end{figure}

\section{Discussion}

Although there are several metric functions to compare the persistence diagrams of large and multiscale data sets, we have found some drawbacks to applying them in social choice data sets. This work is a motivation for a different approach to social choice and welfare science. The most crucial factor on which the comparison of persistence diagrams is based is definitely the accuracy of the persistence points. Their location is relative because of the relative accuracy of invariants' birth and death. The next step to our research is to develop special algorithms for the new field of topological social choice (TSC), which has just emerged.

\par\bigskip\smallskip\par\noindent

\par\noindent
{\it Address}: {\tt {Athanasios Andrikopoulos} }\\ {Department of Computer Engineering \& Informatics\\ University of Patras\\ Greece}
\par\noindent
{\it E-mail address}:{\tt aandriko@ceid.upatras.gr}

\end{document}